\documentclass[a4paper,english]{scrartcl}
\usepackage[utf8x]{inputenc}
\usepackage[T1]{fontenc}
\usepackage{babel}
\usepackage{amsmath,amssymb}
\usepackage[amsmath,thref,thmmarks]{ntheorem}
\usepackage[pdftex]{graphicx}
\usepackage{color}
\usepackage{url}
\usepackage[all]{xy}

\theoremstyle{plain}
\theoremheaderfont{\normalfont\bfseries}
\theorembodyfont{\itshape}
\theoremseparator{:}
\theoremnumbering{arabic}
\theoremsymbol{}

\newtheorem{thm}{Theorem}[section]

\theoremstyle{plain}
\theoremheaderfont{\normalfont\bfseries}
\theorembodyfont{\itshape}
\theoremseparator{:}
\theoremnumbering{arabic}
\theoremsymbol{}

\newtheorem{proposition}[thm]{Proposition}

\theoremstyle{plain}
\theoremheaderfont{\normalfont\bfseries}
\theorembodyfont{\itshape}
\theoremseparator{:}
\theoremnumbering{arabic}
\theoremsymbol{}

\theoremstyle{plain}
\theoremheaderfont{\normalfont\bfseries}
\theorembodyfont{\itshape}
\theoremseparator{:}
\theoremnumbering{arabic}
\theoremsymbol{}

\newtheorem{lemma}[thm]{Lemma}

\theoremstyle{plain}
\theoremheaderfont{\normalfont\bfseries}
\theorembodyfont{\itshape}
\theoremseparator{:}
\theoremnumbering{arabic}
\theoremsymbol{}

\theoremstyle{plain}
\theoremheaderfont{\normalfont\bfseries}
\theorembodyfont{\itshape}
\theoremseparator{:}
\theoremnumbering{Alph}
\theoremsymbol{}

\theoremstyle{plain}
\theoremheaderfont{\normalfont\bfseries}
\theorembodyfont{\itshape}
\theoremseparator{:}
\theoremnumbering{Alph}
\theoremsymbol{}

\theoremstyle{plain}
\theoremheaderfont{\bfseries}
\theorembodyfont{\normalfont}
\theoremseparator{:}
\theoremnumbering{arabic}
\theoremsymbol{}

\newtheorem{definition}[thm]{Definition}

\theoremstyle{plain}
\theoremheaderfont{\bfseries}
\theorembodyfont{\normalfont}
\theoremseparator{:}
\theoremnumbering{arabic}
\theoremsymbol{}

\theoremstyle{plain}
\theoremheaderfont{\bfseries}
\theorembodyfont{\normalfont}
\theoremseparator{:}
\theoremnumbering{arabic}
\theoremsymbol{}

\theoremstyle{nonumberplain}
\theoremheaderfont{\bfseries}
\theorembodyfont{\normalfont}
\theoremseparator{:}
\theoremnumbering{arabic}
\theoremsymbol{}

\newtheorem{bemon}{Remark}

\theoremstyle{nonumberplain}
\theoremheaderfont{\bfseries}
\theorembodyfont{\normalfont}
\theoremseparator{:}
\theoremnumbering{arabic}
\theoremsymbol{}

\theoremstyle{plain}
\theoremheaderfont{\normalfont\scshape}
\theorembodyfont{\itshape}
\theoremseparator{:}
\theoremnumbering{arabic}
\theoremsymbol{}

\newcounter{QuesCounter}
\theoremstyle{plain}
\theoremheaderfont{\normalfont\scshape}
\theorembodyfont{\itshape}
\theoremseparator{:}
\theoremnumbering{arabic}
\theoremsymbol{}

\newtheorem{frage}[QuesCounter]{Question}

\theoremstyle{nonumberplain}
\theoremheaderfont{\normalfont\scshape}
\theorembodyfont{\itshape}
\theoremseparator{:}
\theoremnumbering{arabic}
\theoremsymbol{}

\theoremstyle{nonumberplain}
\theoremseparator{:}
\theoremheaderfont{\scshape}
\theorembodyfont{\upshape}
\theoremsymbol{{\Large\ensuremath{_\dashv}}}

\newtheorem{beweis}{Proof}

\theoremstyle{nonumberplain}
\theoremseparator{:}
\theoremheaderfont{\scshape}
\theorembodyfont{\upshape}
\theoremsymbol{\ensuremath{\square}}

\newcommand{\<}{\langle}
\renewcommand{\>}{\rangle}
\renewcommand{\phi}{\varphi}

\newcommand{\op}[1]{\operatorname{#1}}
\newcommand{\Pot}{\op{\mathcal{P}}}

\newcommand{\GCH}{\op{GCH}}

\newcommand{\AD}{\op{AD}}
\newcommand{\PD}{\op{PD}}

\newcommand{\union}{\bigcup\limits}

\newcommand{\restr}{\upharpoonright}
\newcommand{\concat}{{}^\smallfrown}

\newcommand{\dom}{\op{dom}}

\newcommand{\otp}{\op{otp}}

\newcommand{\ult}{\op{Ult}}
\newcommand{\on}{\op{On}}

\newcommand{\crit}{\op{crit}}
\newcommand{\ran}{\op{ran}}
\newcommand{\id}{\op{id}}
\newcommand{\col}{\op{Col}}
\newcommand{\card}{\op{Card}}
\newcommand{\cof}{\op{cof}}

\newcommand{\Sk}{\op{Sk}}

\newcommand{\fnktsraum}[2]{{}^{#1}{}_{#2}}

\newcommand{\ptwimg}[2]{{#1}"\left[{#2}\right]}
\newcommand{\kleiner}{\mathord{<}}
\newcommand{\kleinergleich}{\mathord{\leq}}

\newcommand{\finsubsets}[1]{{\left[#1\right]}^{\kleiner\omega}}

\newcommand{\eextend}{\trianglelefteq}

\newcommand{\length}{\op{lh}}

\begin{document}

\title{Projective determinacy from long Chang's Conjecture}
\author{Dominik T. Adolf}
\maketitle

\begin{abstract}
 Consider the property $(\aleph_{\omega + 1},\aleph_{\omega + 2},\ldots) \twoheadrightarrow (\aleph_1,\aleph_2,\ldots)$. Here we will show that this property with the addition of the generalized Continuum Hypothesis implies projective determinacy. Of particular interest here is the use of a variant covering argument to prove limited instances of mouse reflection. We believe that this approach could find use for other forms of Chang's Conjecture as well.
\end{abstract}

\section{Introduction}

Let $\lambda,\lambda',\kappa,\kappa'$ be regular cardinals with $\lambda > \kappa$ and $\lambda' > \kappa'$. We take $(\lambda,\kappa) \twoheadrightarrow (\lambda',\kappa')$ to mean: for any structure (in a countable language) on $\lambda$ there exists a substructure $X$ with $\card(X) = \lambda'$ but $\card(X \cap \kappa) = \kappa'$. 

Originating in Model Theory this seemingly innocuous property has significant large cardinal strength. Its most basic form (often known simply as \textit{the} Chang conjecture), $(\aleph_2,\aleph_1) \twoheadrightarrow (\aleph_1,\aleph_0)$ in our notation, is equiconsistent with an $\omega_1$-Erd\H{o}s cardinal.

The exact consistency strength of the analogous property $(\aleph_3,\aleph_2) \twoheadrightarrow (\aleph_2,\aleph_1)$ has proved more elusive. Though it is known to be consistent relative to the existence of a huge cardinal by an argument of Kunen's \cite{kunensatideals}.

Here we will consider what we term a ``long" Chang's conjecture as it involves not two or three but infinitely many cardinals. The degree of strength we can achieve is significantly above what is usually possible for these kind of properties, like in \cite{seancc} or \cite{ralfccone} and \cite{ralfcctwo}. Interestingly enough our long Chang's Conjecture is known to be consistent from a weaker assumption than say $(\aleph_3,\aleph_2) \twoheadrightarrow (\aleph_2,\aleph_1)$, see \cite{ccfromsubcomp}. We believe that some of our methods, specifically those relating to ``mouse reflection", will eventually find use in improving lower bounds for $(\aleph_3,\aleph_2) \twoheadrightarrow (\aleph_2,\aleph_1)$ and other forms of ``short" Chang's Conjecture as well.

$(\aleph_{\omega + 1},\aleph_{\omega + 2},\ldots) \twoheadrightarrow (\aleph_1,\aleph_2,\ldots)$ denotes the property that for any structure (in a countable language) on $\aleph_{\omega + \omega}$  there exists a substructure $X$ with $\card(X \cap \aleph_{\omega + (n + 1)}) = \aleph_{n + 1}$ for any $n < \omega$.

\begin{thm}\label{main}
   Assume $(\aleph_{\omega + 1},\aleph_{\omega + 2},\ldots) \twoheadrightarrow (\aleph_1,\aleph_2,\ldots)$ and $\GCH$. Then projective determinacy holds not only in $V$ but every generic extension of it by a forcing notion of size less than $\aleph_{\omega + \omega}$.
\end{thm}

Projective determinacy is the statement that every game on natural numbers with a projective payoff is determined, \cite{Kechris}. By the foundational results of Descriptive Inner Model theory this statement actually reduces, in fact is equivalent to, the existence of certain canonical inner models. The methods in this paper will stem almost exclusively from the study of such canoncial inner models, ``mice", and only relate to Descriptive Set Theory via these equivalences.

Throughout the paper we will have to assume a general familiarity with Inner Model theory. The purpose of the below subsection is merely to introduce notations and to make clear which particular ``dialect" of inner model theory we shall use throughout this work.

\subsection{Fine structure basics}

We will use the fine structure from \cite{FST}, i.e. extenders are indexed at the local successor of their natural length. $r\Sigma_n$ will be used to denote the levels of the fine structural definability hierarchy. The exact form of the hierarchy depends on the different types of premice as laid out in \cite{FST}, but we shall suppress such detail here. Projecta, notated $\rho_n(\cdot)$, and standard parameters, notated $p_n(\cdot)$, are defined as usual relative to these hierarchies.

 If $\mathcal{M} = (J^{\vec{E}}_{\on \cap \mathcal{M}};\in,\vec{E},F)$ is a premouse we will write $F^\mathcal{M}$ for the top extender, $\mu^\mathcal{M}$ for $\crit(F^\mathcal{M})$ and $\mu^{+,\mathcal{M}}$ for $((\mu^\mathcal{M})^+)^\mathcal{M}$. Note that as a basic feature of the coding $\cof(\mu^{+,\mathcal{M}}) = \cof(\on \cap \mathcal{M})$. For $\alpha \leq (\on \cap \mathcal{M})$ $\mathcal{M} \vert\alpha$ denotes the premouse $(J^{\vec{E}}_\alpha;\in,\vec{E} \restr \alpha,(\vec{E}\concat F)(\alpha))$ while $\mathcal{M} \vert\vert\alpha$ denotes the same structure but with the top predicate removed.
 
Given an iteration tree $\mathcal{T}$ we will write $\mathcal{M}^\mathcal{T}_\beta$ for the $\beta$-th model, $E^\mathcal{T}_\beta$ for the $\beta$-th extender, $\iota^\mathcal{T}_{\alpha,\beta}$ for the embedding from the $\alpha$-th into the $\beta$-th model, and $\mathcal{D}^\mathcal{T}$ for the set of drops. $\leq_\mathcal{T}$ will denote the tree order.
 
Given a tree $\mathcal{T}$ of limit length we will write $\delta(\mathcal{T}):= \sup\limits_{\alpha < \length(\mathcal{T})} \length(E^\mathcal{T}_\alpha)$ and $\mathcal{M}(\mathcal{T}) := \union_{\alpha < \length(\mathcal{T}} \mathcal{M}^\mathcal{T}_\alpha \vert\vert \length(E^\mathcal{T}_\alpha)$. Given a cofinal branch $b$ through $\mathcal{T}$ we will write $\mathcal{Q}(\mathcal{T},b)$ for the least initial segment of $\mathcal{M}^\mathcal{T}_b$ which defines a counterexample to $\delta(\mathcal{T})$ being Woodin or projects below $\delta(\mathcal{T})$ if such exists. By the famous ``zipper" lemma \cite[6.10]{HSTSteel} $\mathcal{Q}(\mathcal{T},b)$ uniquely determines $b$ (for tame mice). This motivates the next definition.

\begin{definition}
   Let $\mathcal{M}$ be a premouse. Let $\mathfrak{M}$ be a mouse operator. An iteration strategy $\Sigma$ is $\mathfrak{M}$-guided if and only if $\mathcal{Q}(\mathcal{T},\Sigma(\mathcal{T})) \eextend \mathfrak{M}(\mathcal{M}(\mathcal{T}))$ for all trees $\mathcal{T} \in \dom(\Sigma)$ of limit length by $\Sigma$.
\end{definition}

Let $n < \omega$. A premouse $\mathcal{M} := (J^{\vec{E}}_\alpha;\in,\vec{E},F)$ is $n$-small if and only if $J^{\vec{E}}_\kappa$ has less than $n$ Woodin cardinals for all $\beta \in \dom(\vec{E}\concat F)$ and $\kappa := \crit((\vec{E}\concat F)_\beta)$.

Let $A$ be a self-wellordered set, i.e. $A$ is wellordered in $L_1(A)$, then $M^\#_n(A)$ denotes the least not $n$-small active premouse over $A$ that is sound above $A$ and every countable substructure is $(\omega_1 + 1)$-iterable.

\begin{lemma}\label{strrefl}
  Let $n < \omega$. Let $\mathcal{M}$ be a premouse of size less than $\kappa$ such that every countable hull has a unique $M^\#_n$-guided $\omega_1$-iteration strategy. Then $\mathcal{M}$ has a unique $M^\#_n$-guided $\kappa$-iteration strategy.
\end{lemma}

\begin{beweis}
  Let $\mathcal{T}$ be an $M^\#_n$-guided iteration tree on $\mathcal{T}$ of size less than $\kappa$. Assume towards a contradiction that $\mathcal{T}$ has no cofinal wellfounded branch $b$ with $\mathcal{Q}(\mathcal{T},b) \eextend M^\#_n(\mathcal{M}(\mathcal{T}))$. Let $H_\theta$ large enough such that it witnesses this fact. Let $X \prec (H_\theta;\in)$ be countable and containing relevant objects. Let $\sigma_X:H_X \rightarrow X$ be an isomorphism where $H_X$ is transitive. By assumption $\mathcal{M}_X := \sigma^{-1}_X(\mathcal{M})$ has a unique $M^\#_n$-guided $\omega_1$-iteration strategy. 
  
  As $\mathcal{T}_X := \sigma^{-1}_X(\mathcal{T})$ is countable there exists a unique cofinal branch $b_X$ with $\mathcal{Q}(\mathcal{T}_X,b_X) \eextend M^\#_n(\mathcal{M}(\mathcal{T}_X))$. Working in $H_X\left[g\right]$ where $g$ is $\col(\omega,\mathcal{T}_X)$-generic over $H_X$, we have that the existence of a branch such as $b_X$ is $\Sigma^1_1$ in real codes for $\mathcal{T}_X$ and $M^\#_n(\mathcal{M}(\mathcal{T}_X))$. By absoluteness we have $b_X \in H_X\left[g\right]$. But $b_X$ is independent of the generic so actually $b_X \in H_X$. Contradiction!
\end{beweis}

Results such as these demonstrate the importance of Woodin cardinals and canonical inner models with Woodin cardinals for the study of iterability. But they also have a profound impact on the determinacy of certain games on natural numbers.

\begin{thm}[Woodin, Martin-Steel]
    Projective Determinacy holds if and only if $M^\#_n(x)$ exists, i.e. is $\omega_1$-iterable, for all reals $x$.
\end{thm}

So, in order to prove the main theorem it will be sufficient to prove the existence of $M^\#_n(x)$ for all reals $x$ from the assumption of $(\aleph_{\omega + 1},\aleph_{\omega + 2},\ldots) \twoheadrightarrow (\aleph_1,\aleph_2,\ldots)$ and the $\GCH$. In fact, we will prove the existence of such for all bounded subsets of $\aleph_{\omega + \omega}$ which shows the generic absoluteness of $\PD$ for small forcing notions as $M^\#_n$'s determine themselves on generic extensions, see \cite[1.4.10]{CMI}.

The pursuit of this goal which will take up the entirety of this paper will be divided in the following parts: we begin with a short introduction (you are here); we will discuss some basic facts about special substructures,the existence of which will be guaranteed by our assumption; we will then first, as a warmup, show the existence of an inner model with a Woodin cardinal and proceed from there to the full Core Model Induction argument; most interesting part of which will be in how to achieve mouse reflection; we do not see how to make an abstract argument work (say as in \cite{PFAimplADLR}), instead we will appeal to the Core Model Dichotomy Theorem to show mouse reflection for $M^\#_{n + 1}$'s; this will require an alternative covering argument in which $M^\#_{n + 1}$ takes the role of the core model, $K$.

\subsection{Acknowlegdements}

The author would like to acknowledge the support of Omer Ben-Neria and his grant from the Israel Science Foundation (Grant 1832/19).

\section{Long Chang's Conjecture basic facts}

Generally, we will use the language of stationary sets throughout this paper. Recall, a subset $S$ of $\Pot(\mathcal{X})$ is stationary iff for all $F:\finsubsets{\mathcal{X}} \rightarrow \mathcal{X}$ there is some $X \subseteq \mathcal{X}$ in $S$ that is closed under $F$. $(\aleph_{\omega + 1},\aleph_{\omega + 2},\ldots) \twoheadrightarrow (\aleph_1,\aleph_2,\ldots)$ then means: there is a stationary set of $X \subset \aleph_{\omega + \omega}$ with $\card(X \cap \aleph_{\omega + (n + 1)}) = \aleph_{n + 1}$ for all $n < \omega$. 

Note that this is equivalent to saying that there is a stationary set of $X \subset H_{\aleph_{\omega + \omega}}$ with $X \cap \aleph_{\omega + \omega}$ satisfying $(\aleph_{\omega + 1},\aleph_{\omega + 2},\ldots) \twoheadrightarrow (\aleph_1,\aleph_2,\ldots)$.

For the arguments to come it would be a great help if our hulls were locally countably closed. This requirement can likely be weakend slightly. Unfortunately, we do not know if it can be removed completely as in \cite{OpistolfromJonsson}.

\begin{definition}
    A structure $X \prec (H_{\omega + \omega};\in)$ is appropriate if and only if $X$ witnesses $(\aleph_{\omega + 1},\aleph_{\omega +  2},\ldots) \twoheadrightarrow (\aleph_{\omega_1},\aleph_{\omega_2}.\ldots)$, i.e. $\card(X \cap \aleph_{\omega + (n + 1)}) = \aleph_{n + 1}$ for all $n < \omega$, and  $\fnktsraum{\omega}{X \cap \aleph_\omega} \subset X$.
\end{definition}

\begin{lemma}
	Assume $(\aleph_{\omega + 1},\aleph_{\omega + 2},\ldots) \twoheadrightarrow (\aleph_1,\aleph_2,\ldots)$ and the $\GCH$. Then there exists a stationary set of appropriate structures.
\end{lemma}

\begin{beweis}
	Let $\mathfrak{A}$ be a Skolemized structure on $H_{\aleph_{\omega + \omega}}$. By assumption, we know that there exists some $X \prec \mathfrak{A}$ with $\card(X \cap \aleph_{\omega + (n + 1)}) = \aleph_{n + 1}$ for all $n < \omega$. We shall show that for all $x \subset \aleph_\omega$ of size $\kleinergleich \aleph_1$ the Skolemhull of $X \cup x$ (notated here as $\Sk^{\mathfrak{A}}(X \cup x)$) is a Chang structure of the required type. By iterating this procedure we can then construct an $\subset$-increasing and continuous sequence of Chang structures $\<X_\gamma:\gamma \leq \omega_1\>$ with $X_\gamma \prec \mathfrak{A}$ and $\fnktsraum{\omega}{X_\beta \cap \aleph_\omega} \subset X_\gamma$ for all $\beta < \gamma \leq \aleph_1$. 
	
	$X_{\aleph_1}$ then is a substructure of $\mathfrak{A}$ of the required type with $\fnktsraum{\omega}{X \cap \aleph_\omega} \subset X$. 
	
	Consider then $Y := \Sk^\mathfrak{A}(X \cup x)$. Let $n < \omega$. Without loss of generality, we can assume that $\mathfrak{A}$ satisfies the collection scheme. Therefore we have that any element of $Y \cap \aleph_{\omega + (n + 1)}$ can be written as $f(\beta)$ where $f:\aleph_\omega \rightarrow \aleph_{\omega + (n + 1)}$,$f \in X$ and $\beta \in x$.
	
	By our assumption there are only $\aleph_{\omega + (n + 1)}$ many functions from $\aleph_\omega$ into $\aleph_{\omega + (n + 1)}$. We then have some bijective $g:\aleph_{\omega + (n + 1)} \rightarrow \fnktsraum{\aleph_\omega}{\aleph_{\omega + (n + 1)}}$ in $H_{\aleph_{\omega + \omega}}$. By elementarity, there is some such $g \in X$ witnessing that $\card(\fnktsraum{\aleph_\omega}{\aleph_{\omega + (n + 1)}} \cap X) = \card(X \cap \aleph_{\omega + (n + 1)}) = \aleph_{n + 1}$. We conclude then that $\card(Y \cap \aleph_{\omega + (n + 1)}) = \aleph_{n + 1} \cdot \aleph_1 = \aleph_{n + 1}$ as desired.
\end{beweis}

\begin{lemma}
    $\otp(X \cap \aleph_{\omega + (n + 2)}) = \aleph_{n + 2}$ for all $n < \omega$, for all but non-stationarily many appropriate structures $X$.
\end{lemma}

\begin{beweis}
    Let $n < \omega$. Let $f: (\aleph_{\omega + \omega})^2 \rightarrow (\aleph_{\omega + \omega})$ be a function such that $f \restr \{\alpha\} \times \alpha$ is an injection into $\card(\alpha)$ for all $\alpha < \aleph_{\omega + \omega}$. Any appropriate $X$ closed under $f$ will do. Assume $\otp(X \cap \aleph_{\omega + (n + 2)}) > \aleph_{n + 2}$, let $\alpha < \aleph_{\omega + (n + 2)}$ be such that $\otp(X \cap \alpha) = \aleph_{n + 2}$. $f \restr \{\alpha\} \times (X \cap \alpha)$ is then an injection into $X \cap \aleph_{\omega + (n + 1)}$. Therefore $\card(X \cap \aleph_{\omega + (n + 1)}) \geq \aleph_{n + 2}$. Contradiction!
\end{beweis}

Important here is mainly the implication that $\cof(X \cap \aleph_{\omega + (n + 2)}) = \aleph_{n + 2}$. The same actually holds for $n+1$. The argument here goes back to a theorem of Shelah's about cofinalities of successor cardinals in generic extensions \cite[p 451]{Jech}. Foreman and Magidor noticed that a similar property holds for substructures \cite[2.15]{fandmdefcounterexampletoch}.

\begin{lemma}
   $\cof(X \cap \aleph_{\omega + 1}) = \omega_1$ for all but non-stationarily many appropriate structures $X$.
\end{lemma}

\begin{beweis}
     Let $\vec{x} := \<x_\alpha: \alpha < \aleph_{\omega + 2}\>$ be an almost disjoint family of subsets of $\aleph_{\omega + 1}$. Assume $\vec{x} \in X$. Let $g:\omega_1 \rightarrow (X \cap \aleph_{\omega + 1})$ be a bijection. Assume for a contradiction that $\cof(X \cap \aleph_{\omega + 1}) = \omega$. Then for every $\alpha \in X \cap \aleph_{\omega + 2}$ there exists some $\beta_\alpha < \omega_1$ such that $\ptwimg{g}{\beta_\alpha}$ is unbounded in $x_\alpha \cap X$. There must be an unbounded $Y \subset X \cap \aleph_{\omega + 2}$ and some fixed $\beta$ such that $\beta_\alpha = \beta$ for all $\alpha \in Y$. Let $\gamma \in Y$ be such that $\otp(Y \cap \gamma) = \omega_1$.

Note then that there exists a sequence $\vec{\delta} := \<\delta_\alpha: \alpha < \gamma\>$ such that $\<x_\alpha \backslash \delta_\alpha: \alpha < \gamma\>$ is pairwise disjoint. (To see this: reorganize $\<x_\alpha : \alpha < \gamma\>$ in ordertype $\aleph_{\omega + 1}$ and use the regularity of $\aleph_{\omega + 1}$.) We can assume that $\vec{\delta} \in X$. Crucially, this implies $\delta_\alpha \in X \cap \aleph_{\omega + 1}$ whenever $\alpha \in X$. We then have that $\<\ptwimg{g^{-1}}{x_\alpha \backslash \delta_\alpha} \cap \beta : \alpha \in Y \cap \gamma\>$ is a pairwise disjoint family of size $\omega_1$ that consists of subsets of $\beta$. Contradiction!
\end{beweis}

\begin{definition}
   An appropriate structure $X \prec (H_{\aleph_{\omega + \omega}};\in)$ is good if and only if $\cof(X \cap \aleph_{\omega + (n + 1)}) = \aleph_{n + 1}$ for all $n < \omega$.
\end{definition}

\begin{lemma}\label{closurelemma}
   Let $X$ be good. Then $\fnktsraum{\omega}{X \cap \aleph_{\omega + n}} \subset X$ for all $n < \omega$.
\end{lemma}

\begin{beweis}
    Proof is by induction. $n = 0$ is by assumption. Let then $x \subset (X \cap \aleph_{\omega + {n + 1}})$ be countable. As $X$ is good there is some $\alpha \in (X \cap \aleph_{\omega + (n + 1)})$ such that $x \subset \alpha$. Pick some bijection $g: \aleph_{\omega + n} \rightarrow \beta$ in $X$. By induction assumption we have $y := \ptwimg{g^{-1}}{x}$ is a countable subset of $X \cap \aleph_{\omega + n}$ and hence in $X$. Then $x = \ptwimg{g}{y}$ is in $X$.
\end{beweis}

We will close this section with a short discussion on square principles which will have some minor use in section 4.

\begin{definition}
   Let $\kappa$ be a cardinal. The principle $\square_\kappa$ holds if there exists a sequence $\<C_\alpha: \alpha \in \kappa^+ \cap \lim\>$ such that
   \begin{enumerate}
      \item $C_\alpha \subset \alpha$ is a club for all $\alpha \in \kappa^+ \cap \lim$,
      \item $\beta \in \lim(C_\alpha) \Rightarrow C_\beta = C_\alpha \cap \beta$ for all $\beta < \alpha \in \kappa^+ \cap \lim$,
      \item $\otp(C_\alpha) \leq \kappa$ for all $\alpha \in \kappa^+ \cap \lim$.
   \end{enumerate}
\end{definition}

\begin{lemma}\label{nosquarelemma}
   Assume $(\aleph_{\omega + 1},\aleph_{\omega + 2},\ldots) \twoheadrightarrow (\aleph_1,\aleph_2,\ldots)$ and the $\GCH$. Then $\square_{\aleph_{\omega + (n + 1)}}$ fails for all $n < \omega$.
\end{lemma}

\begin{beweis}
    Let $n < \omega$. Assume for a contradiction that $\vec{C} := \<C_\alpha : \alpha \in \aleph_{\omega + (n + 2)} \cap \lim\>$ is a $\square_{\aleph_{\omega + (n + 1)}}$-sequence. Let $X \prec (H_{\aleph_{\omega + \omega}};\in)$ be a good structure with $\vec{C} \in X$ and $\otp(X \cap \aleph_{\omega + (n + 2)}) = \aleph_{n + 2}$. Let $\sigma_X: H_X \rightarrow X$ be an isomorphism where $H_X$ is transitive. Let $\vec{D} = \sigma^{-1}_X(\vec{C})$. Note that $\sigma_X$ is continuous at points of countable cofinality. Hence $\ran(\sigma_X) \cap \aleph_{\omega + (n + 2)}$ is an $\omega$-club. Therefore $D : = \ptwimg{\sigma^{-1}_X}{C_{\sup(X \cap \aleph_{\omega + (n + 2)})}}$ is a set of ordertype $\aleph_{n + 2}$.

We are finished if we can show that $D_\alpha = D \cap \alpha$ for unboundedly many $\alpha < \aleph_{n + 2}$ as this implies that $D_\alpha$ have ordertypes unbounded in $\aleph_{n + 2}$ which contradicts the fact that by elementarity ordertypes are bounded by $\sigma^{-1}_X(\aleph_{\omega + (n + 1)}) < \aleph_{n + 2}$. Take then some arbitrary $\alpha$ which is a limit point of $D$ of countable cofinality. Then by continuity $\sigma_X(\alpha)$ is a limit point of $C_{\sup(X \cap \aleph_{\omega + (n + 2)})}$ and hence $C_{\sigma_X(\alpha)} = C_{\sup(X \cap \aleph_{\omega + (n + 2)})} \cap \sigma_X(\alpha)$. This implies
\[ D \cap \alpha = \ptwimg{\sigma^{-1}_X}{C_{\sup(X \cap \aleph_{\omega + (n + 2)})} \cap \sigma_X(\alpha)} = \ptwimg{\sigma^{-1}_X}{C_{\sigma_X(\alpha)}} = D_\alpha. \]
Contradiction!
\end{beweis}

We do not see if something similar might apply to the principle $\square(\aleph_{\omega + (n + 2)})$. If we could this would allow us to point to \cite{cohseqthread} to get significant strength immediately. Another question is if $\square_{\aleph_{\omega + \omega}}$ must fail under the same assumption which also would give strength immediately by \cite{PFAimplADLR}. But this fails.

\begin{lemma}
   Assume $(\aleph_{\omega + 1},\aleph_{\omega + 2},\ldots) \twoheadrightarrow (\aleph_1,\aleph_2,\ldots)$ and the $\GCH$. Then there exists an inner model in which both of these properties and additionally $\square_{\aleph_{\omega + \omega}}$ holds.
\end{lemma}

\begin{beweis}
    The inner model is $L(H_{\aleph_{\omega + \omega}})$. By the $\GCH$ the set $H_{\aleph_{\omega + \omega}}$ is coded by a subset of $\aleph_{\omega + \omega}$ and hence by Jensen's proof $\square_{\aleph_{\omega + \omega}}$ holds in $L(H_{\aleph_{\omega + \omega}})$, and so does the $\GCH$, see \cite{HSTRalfMartin}. It remains to be seen that $(\aleph_{\omega + 1},\aleph_{\omega + 2},\ldots) \twoheadrightarrow (\aleph_1,\aleph_2,\ldots)$ holds there. Fix some $F:\finsubsets{\aleph_{\omega + \omega}} \rightarrow \aleph_{\omega + \omega}$ in $L(H_{\aleph_{\omega + \omega}})$. We build a tree of finite sequences $\<X_0,\ldots,X_{n - 1}\>$ such that

   \begin{itemize}
     \item $X_k \prec (H_{\aleph_{\omega + (k + 1)}};\in)$ for all $k < n$,
     \item $X_k = X_l \cap H_{\aleph_{\omega + (k + 1)}}$ for all $k < l < n$,
     \item $\card(X_k) = \aleph_{k + 1}$ for all $k < n$,
     \item $F(a) < \aleph_{\omega + (k + 1)} \Rightarrow F(a) \in X_k$ for all $a \in \finsubsets{X \cap \aleph_{\omega + (k + 1)}}$, all $k < n$.
   \end{itemize}

   By $(\aleph_{\omega + 1},\aleph_{\omega + 2},\ldots) \twoheadrightarrow (\aleph_1,\aleph_2,\ldots)$ in $V$ this tree is ill-founded. This also holds in $L(H_{\aleph_{\omega + \omega}})$. Any witness to this ill-foundedness will witness $(\aleph_{\omega + 1},\aleph_{\omega + 2},\ldots) \twoheadrightarrow (\aleph_1,\aleph_2,\ldots)$ and be $F$-closed.
\end{beweis}

\section{Towards one Woodin cardinal}

For this section assume that there is no inner model with a Woodin cardinal. This implies that there exists a Jensen-Steel core model \cite{knombl}. Let $W$ be a soundness witness for $K$ up to $\aleph_{\omega + \omega}$. $W$ is a weasel of height $\kappa^+$ where $\kappa$ is a singular strong limit of cofinality above $\aleph_{\omega \cdot 2 + 1}$. $W$ has the definability property (relative to $\aleph_{\omega \cdot 2 + 1}$-club subsets of $\kappa^+$) at all points below $\aleph_{\omega + \omega}$.

Our goal for this section will be to show that $(\aleph_{\omega + 1},\aleph_{\omega + 2},\ldots) \twoheadrightarrow (\aleph_1,\aleph_2,\ldots)$ and the $\GCH$ cannot hold together under this assumption. Towards a contradiction assume then that they do hold. Let then $X \prec (H_{\kappa^{++}};\in)$ be a structure such that $W \in X$ and $X \cap H_{\omega + \omega}$ is good. Let $\sigma_X:H_X \rightarrow X$ be an isomorphism where $H_X$ is transitive. Let $W_X := \sigma^{-1}_X(W)$. We will co-iterate $W_X$ and $W$ but only until they are lined up to $\aleph_\omega$. Let $W^*_X$ be the final model on the $W_X$ side of the iteration (we will eventually realize that $W^*_X = W_X$). We let $\mathcal{T}$ be the tree on the $W$-side of the iteration.

We will write $\mathcal{M}_\beta$ for the $\beta$-th model of $\mathcal{T}$, $E_\beta$ for the $\beta$-th extender, $\iota_{\beta,\gamma}: \mathcal{M}_\beta \rightarrow \mathcal{M}_\gamma$ for the iteration embedding, and $\mathcal{D} \subset \length(\mathcal{T})$ for the set of drops. We would usually write $\mathcal{M}^\mathcal{T}_\beta$ etc, but $\mathcal{T}$ will be the only tree this section so no confusion should arise.

The first main consequence of goodness will be the following effect on the cofinalites of successor cardinals in $W_X$:

\begin{lemma}\label{sec3cofincore}
   Let $\kappa < \aleph_\omega$ be above $\crit(\sigma_X)$. Then $\cof((\kappa^+)^{W_X} = \card(\kappa)$.
\end{lemma}

\begin{beweis}
    By assumption $\sigma_X(\kappa)$ is above $\aleph_2$ and below some $\aleph_{\omega + (n + 2)}$ for some $n$. By \thref{nosquarelemma} and \cite{square} the latter are limit cardinals in $W$, thus $\card(\sigma_X(\kappa)) = \card((\sigma_X(\kappa)^+)^W)$. By weak covering we have $\cof(\sigma_X((\kappa^+)^{W_X}) = \card(\sigma_X(\kappa))$. Some witness to this can be assumed to be an element of $X$. We then have $\cof(\sigma_X((\kappa^+)^{W_X} \cap X) = \cof(X \cap \card(\sigma_X(\kappa)))$.  The latter equals $\card(X \cap \card(\sigma_X(\kappa)))$ by goodness.
\end{beweis}

The preceding lemma used the weak covering property in $W$. To continue we will have to consider the proof of this property from \cite{coveringbelowwoodin}. Thus we will consider three phalanxes: one expanding $\mathcal{T}$ and two derived from the first one. The index set of those phalanxes shall be the set of infinite cardinals of $W^*_X$, but limited to those below $\aleph_\omega$. For an ordinal $\alpha < \aleph_\omega$, let $\kappa_\alpha := \aleph^{W^*_X}_\alpha$ and $\lambda_\alpha := \aleph^{W^*_X}_{\alpha + 1}$.

We will now describe the objects making up the aforementioned three phalanxes:

\begin{itemize}
    \item $\mathcal{P}_\alpha$ is the least level of $\mathcal{M}_\beta$ that agrees with $W^*_X$ up to $\lambda_\alpha$ and projects to $\kappa_\alpha$ or below (if such exists, otherwise $\mathcal{P}_\alpha = \mathcal{M}_\beta$) where $\beta < \length(\mathcal{T})$ is least with $\nu(E_\beta) > \kappa_\alpha$ (if such exists, otherwise $\beta = \length(\mathcal{T}) - 1$);
    \item $n_\alpha$ is least such that $\rho^{\mathcal{P}_\alpha}_{n_\alpha + 1} \leq \kappa_\alpha$ (if such exists, otherwise $n_\alpha = \omega$);
    \item $\mathcal{R}_\alpha := \ult_{n_\alpha}(\mathcal{P}_\alpha;\sigma_X \restr W_X \vert\vert \lambda_\alpha)$, let $\sigma_\alpha$ be the ultrapower embedding;
    \item if $\mathcal{P}_\alpha$ is an active type I or II premouse such that $(\mu^+)^{\mathcal{P}_\alpha} \leq \kappa_\alpha$ is a discontinuity point of $\sigma_X$ and $n_\alpha = 0$, then $\mathcal{R}_\alpha$ is a proto-mouse; in that case define $\mathcal{S}_\alpha := \ult_{k}(\mathcal{S}_\beta;F^{\mathcal{R}_\alpha})$ where $\beta,k$ are such that
          \begin{itemize}
              \item there is a descending sequence $\gamma_0,\ldots,\gamma_m$ such that $\gamma_0 = \alpha$ and $\gamma_{n+1}$ is the unique ordinal $\gamma$ with $\kappa_\gamma = \mu^{ \mathcal{P}_{\gamma_n}}$,
              \item $m$ is the first number such that $\mathcal{R}_{\gamma^X_m}$ is a premouse,
              \item $\beta = \gamma_1$ and $k = n_{\gamma_m}$,
          \end{itemize}
    otherwise $\mathcal{S}_\alpha := \mathcal{R}_\alpha$.
\end{itemize}

There is also the mouse $\mathcal{Q}_\alpha$ which is defined from $\mathcal{P}_\alpha$ the same way that $\mathcal{S}_\alpha$ is defined from $\mathcal{R}_\alpha$. Its crucial property is that $\mathcal{S}_\alpha$ can be realized as a long ultrapower of it. It will not be relevant for this section but we will have to go into more detail for the next.

In \cite{coveringbelowwoodin} six properties $(1)_\alpha$ up to $(6)_\alpha$ of these objects are laid out. They are:

\begin{itemize}
    \item $(1)_\alpha$ is the statement that the least extender applied to $W_X$ in the co-iteration of $W_X$ and $W$ has length at least $\lambda_\alpha$;
    \item $(2)_\alpha$ is the statement that the phalanx $\<\<W,\mathcal{S}_\alpha\>,\sigma_X(\kappa_\alpha)\>$ is iterable;
    \item $(3)_\alpha$ is the statement that the phalanx $\<\<W_X,\mathcal{Q}_\alpha\>,\kappa^X_\alpha\>$ is iterable;
    \item $(4)_\alpha$ is the statement that the phalanx $\<\<\mathcal{P}_\beta: \beta < \alpha\> \concat W_X\>,\<\lambda_\beta: \beta < \alpha\>\>$ is iterable;
    \item $(5)_\alpha$ is the statement that the special phalanx $\<\<\mathcal{R}_\beta: \beta < \alpha\> \concat W\>,\<\chi_X(\sigma_X(\lambda_\beta)):\beta < \alpha\>\>$ is iterable;
    \item $(6)_\alpha$ is the statement that the special phalanx $\<\<\mathcal{S}_\beta: \beta < \alpha\> \concat W\>,\<\chi_X(\sigma_X(\lambda_\beta)):\beta < \alpha\>\>$ is iterable.
\end{itemize}

\begin{bemon}
    Note that the phalanxes from $(5)_\alpha$ and $(6)_\alpha$ are ``special", i.e. the exchange ordinals are not necessarily cardinals in subsequent models. This obviously necessitates a great deal of care in their use. Generally, these phalanxes are not intended to be used in comparisons so we shall not run into trouble using them.
\end{bemon}

From these we really only need $(1)_\alpha$ for all $\alpha$ such that $\kappa_\alpha < \aleph_\omega$, but the only way to get there seems by an interlaced induction using all these properties which is why, especially in the next section, we will have to consider all of them. We will also need a method to show that certain of the $\mathcal{P}_\alpha$ are not weasels for which we will need $(2)_\alpha$. We could avoid these by using $\lambda$-indexing but our sources use ms-indexing so we do the same.

By \cite{coveringbelowwoodin} these properties satisfy the following relationships:

\begin{proposition}[Mitchell-Schimmerling-Steel]\label{coveringinduction}
Assume that $\alpha$ is such that $(1)^X_\beta$ holds for all $\beta < \alpha$, then:
    \begin{itemize}
        \item[$(i)$] $(6)_\alpha \Rightarrow (5)_\alpha \Rightarrow (4)_\alpha \Rightarrow (1)_\alpha$;
        \item[$(ii)$] $\forall \beta < \alpha (4)_\beta \Rightarrow (3)_\alpha$, $\forall \beta < \alpha (2)_\beta \rightarrow (6)_\alpha$.
    \end{itemize}
\end{proposition}

So to have all properties everywhere we need to have that $(3)_\alpha$ implies $(2)_\alpha$. This is where the authors of \cite{coveringbelowwoodin} use the countable closure of their hull. We have slightly less than that. Recall, goodness implies that $\fnktsraum{\omega}{\aleph_{\omega + (n + 1)} \cap X} \subset X$ for all $n < \omega$. This will be enough to show this for $\alpha$ such that $\kappa_\alpha < \aleph_\omega$. 

\begin{lemma}
     Let $\alpha$ be such that $\kappa_\alpha < \aleph_\omega$ and $(1)_\beta$ holds for all $\beta < \alpha$. Then $(3)_\alpha$ implies $(2)_\alpha$.
\end{lemma}

\begin{beweis}
   Assume $(2)_\alpha$ fails. Take then a countable hull $Y$ of the phalanx witnessing this failure. $Y \cap \mathcal{S}_\alpha$ can be represented by a sequence $\<\left[g_n,b_n\right]_{\sigma_X \restr \lambda_\alpha}: n < \omega\>$. The relations between these elements can thus be coded by sets
   \[ A_{\phi,a} := \sigma_X(\{c \in \left[\sigma_X(\kappa_\alpha)\right]^{\vert d\vert} \vert \mathcal{Q}_\alpha \models \phi(g_n(c^{d,b_n}): n \in a)\})\]
   where $\phi$ is a first-order formula in the language of premice, $a \in \finsubsets{\omega}$ and $d = \bigcup\limits_{n \in a} b_n$. Each of these sets is an element of $W \vert\vert \chi_X(\sigma_X(\lambda_\alpha))$ and is thus represented by a finite set of ordinals. Our local closure is thus sufficient to show that the sequence of all such sets is in $X$. As usual this can be used to show that $Y \cap \mathcal{S}_\alpha$ can be realized into $\mathcal{Q}_\alpha$.
   
   The same trick does not work for $Y \cap W$ as the weasel is larger than our degree of closure. Instead take some $f \in X$ with $f: \Omega \twoheadrightarrow W$. Then $Y \cap W$ can be represented by a sequence $\<f(a_n): n < \omega\>$. The relations between these elements can thus be coded by sets
   \[ B_{\phi,a} := \{c \in \left[\Omega\right]^{\vert d\vert}\vert W \models \phi(f(c^{d,a_n}): n \in a)\}\]
   where $\phi$ is a first order formula in the language of premice, $a \in \finsubsets{\omega}$ and $d = \bigcup\limits_{n \in a} a_n$. Note how $B_{\phi,a}$ only depends on $f$ and a finite type. Thus the sequence of these sets is definable from $f$ and a sequence of such types, and is thus an element of $X$.
   
   Thus we can also realize $Y \cap W$ into $W_X$. If we could show that the two realizations can be chosen to have sufficient agreement, then we reach a contradiction. To do that $X$ needs to have a little bit more information. We need to ensure that if $f(a_n)$ and $g_m(b_m)$ represent the same ordinal below $\chi_X(\sigma_X(\lambda_\alpha))$ then the same is true for our realizations. For this we need the set
   \[ C_{n.m} = \{c \in \left[\Omega\right]^{\vert d\vert}\vert f(c^{d,a_n}) = \sigma_X(g_m)(c^{d,b_m})\}\]
   where $d = a_n \cup b_m$. Crucially, we can in this special case choose $g_m$ such that $\sigma_X(g_m) \in X$. Note that this set is definable from $f,\sigma_X(g_m)$ and a finite type. As $\sigma_X(g_m)$ is coded by a finite set of ordinals below $\sigma_X(\lambda_\alpha)$ we have enough closure to show that the sequence is in $X$. 
\end{beweis}

The importance of the property $(2)_\alpha$ for us in particular is:

\begin{proposition}[Schindler]\label{collapsinglevelsII}
   Assume $\alpha$ be such that $(2)_\alpha$ holds. If $\on \cap \mathcal{S}_\alpha = \Omega$ then there exists an extender $G_\alpha$ on the $W$-sequence such that $S_\alpha = \ult(W;G^X_\alpha)$, additionally, $( \length(G^X_\alpha),0)$ is the lexicographically least $(\gamma,m)$ such that a surjection from $\sigma_X(\kappa_\alpha)$ onto $\chi_X(\sigma_X(\lambda_\alpha))$ is $r\Sigma_{m + 1}$-definable over $W\vert \gamma$.
\end{proposition}

This is actually a strengthening of what was proved in \cite{coveringbelowwoodin}. See \cite{gitikshelahschindlerpcfwoodincard}.

We can now attempt the main argument of this section: let $\alpha_n$ be the unique $\alpha$ such that $\kappa_\alpha = \aleph_{n + 2}$, let $\xi_n$ be the least $\xi$ such that $\nu(E_\xi) \geq \aleph_{n+2}$ (if such exists, otherwise $\xi_n := \length(\mathcal{T}) - 1$).

\begin{lemma}\label{sec3cutpoint}
    $\beta \geq \xi_n \Rightarrow \beta \geq_\mathcal{T} \xi_n$ for all $n < \omega$.
\end{lemma}

\begin{beweis}
    Assume not. Then there exists some $\beta \geq \xi_n$ with $\length(E_\beta) \geq \aleph_{n+2}$ but $\crit(E_\beta) < \aleph_{n + 2}$. As $W_X$ does not move in the iteration $\length(E_\beta)$ is a successor cardinal in $W_X$ above $\aleph_{n+2}$. By \thref{sec3cofincore} $\cof(\length(E_\beta)) \geq \aleph_{n+2}$. The former equals $\cof((\crit(E_\beta)^+)^{W_X})$, but by \thref{nosquarelemma} and \cite{square} $\aleph_{n+2}$ is a limit cardinal in $W_X$. Hence $(\crit(E_\beta)^+)^{W_X} < \aleph_{n + 2}$. Contradiction!
\end{beweis}

So $\{\xi_n : n < \omega\} \subset \left[0,\length(\mathcal{T}) - 1\right)$. Next we will see that there must be at least one drop on this branch.

\begin{lemma}\label{sec3firstdrop}
    $\mathcal{P}_{\alpha_n}$ is not a weasel for all $n < \omega$.
\end{lemma}

\begin{beweis}
    Assume not. Then $\mathcal{S}_{\alpha_n} = \mathcal{R}_{\alpha_n}$ is also a weasel. Thus we can apply \thref{collapsinglevelsII} and get some extender $G_{\alpha_n}$ such that $\mathcal{S}_{\alpha_n} = \ult(W;G_{\alpha_n})$. A crucial thing to realize here is that $(\aleph^+_{n+2})^{\mathcal{S}_{\alpha_n}} = \chi_X((\aleph^+_{\omega + (n + 2)})^W)$. By \thref{collapsinglevelsII} then $\cof(\length(G_{\alpha_n})) = \cof((\aleph^+_{n + 2})^{\mathcal{S}_{\alpha_n}})$. The latter equals $\aleph_{n+ 2}$ by \thref{sec3cofincore}, but the former equals $\cof((\crit(G_{\alpha_n})^+)^W)$.

Note that there are two embeddings $W \rightarrow \mathcal{S}_{\alpha_n}$: $\iota_{G_{\alpha_n}}$, and $\sigma_{\alpha_n} \circ \iota_{0,\xi_n}$. By the definability property these embeddings are equal. We conclude that $\crit(G_{\alpha_n}) < \aleph_2$. Moreover, its successor too must be below $\aleph_2$ as this is a limit in $W_X$. Contradiction!
\end{beweis}

Note then that some initial segment of $\mathcal{M}_{\xi_0}$ projects below $\aleph_2$. On the other hand any extender applied to $\mathcal{M}_{\xi_0}$ in the course of the iteration has critical point at least $\aleph_2$. Hence $\mathcal{D} \cap \left[0,\xi_1\right)_\mathcal{T} \neq \emptyset$. This will be enough to show that there are in fact infinitely many more drops on the main branch.

\begin{lemma}\label{sec3moredrops}
   $\mathcal{D} \cap \left[\xi_{n + 1},\xi_{n + 2}\right)_\mathcal{T} \neq \emptyset$ for all $n < \omega$.
\end{lemma}

\begin{beweis}
   By previously established facts $\mathcal{M}_{\xi_{n + 1}}$ is the direct limit of mice that are of size less than $\aleph_{n + 3}$. Therefore $\cof(\rho_n(\mathcal{M}_{\xi_{n+1}})) < \aleph_{n+3}$ where $n$ is the degree of the branch $\left[0,\xi_{n + 1}\right)_\mathcal{T}$. As $\mathcal{M}_{\xi_{n+1}}$ is $(n+1)$-sound above $\aleph_{n+3}$ we have $\cof((\aleph^+_{n+3})^{\mathcal{M}_{\xi_{n+1}}}) < \aleph_{n + 3}$. But $\cof((\aleph^+_{n+3})^{W_X}) = \aleph_{n + 3}$ by \thref{sec3cofincore}. As no extenders get applied on the $W_X$-side, we must have $(\aleph^+_{n+3})^{W_X} < (\aleph^+_{n+3})^{\mathcal{M}_{\xi_n}}$. By \thref{sec3cutpoint} then any node coming out of $\xi_{n+1}$ must be a drop.
\end{beweis}

So, in conlusion, we have infinitely many drops occuring on the final branch of $\mathcal{T}$. But this contradicts the iterabilty of $W$!

\section{Towards projective determinacy}

Our goal for this section will be to show that $M^\#_n(A)$ exists for all bounded subsets $A \subset \aleph_{\omega + \omega}$, assuming that $(\aleph_{\omega + 1},\aleph_{\omega + 2},\ldots) \twoheadrightarrow (\aleph_1,\aleph_2,\ldots)$ and the $\GCH$ holds.

We must first start with regular sharps then. Here the arguments is also significantly easier. Assume then for the remainder of the section that $(\aleph_{\omega + 1},\aleph_{\omega + 2},\ldots) \twoheadrightarrow (\aleph_1,\aleph_2,\ldots)$ and the $\GCH$ holds.

\begin{lemma}
   Let $A \subset \aleph_1$. Then $A^\#$ exists.
\end{lemma}

\begin{beweis}
   Let $X \prec (H_\theta,\in)$ be such that $\theta >> \aleph_{\omega+ \omega}$, $A \in X$, and $X \cap H_{\omega + \omega}$ is good. Let $\sigma_X:H_X \rightarrow X$ be an isomorphism where $H_X$ is transitive. By condensation $(L(A))^{H_X} = L_{\on \cap H_X}(A)$. We then have $\sigma_X \restr (L(A))^{H_X} \rightarrow L_\theta(A)$. As $\aleph_2 \subset H_X$ this is enough to get $A^\#$.
\end{beweis}

The rest of the argument from here then ammounts to proving a particular instance of mouse reflection. Unlike in other arguments, e.g. \cite{PFAimplADLR}, we do not seem to have access to an abstract argument that works for all mouse operators. Instead the argument must be adapted to each individual mouse operator.

\begin{lemma}
     Let $A \subset \aleph_{\omega + \omega}$ be bounded. If $\cof((\card(A)^+)^{L(A)}) = \omega$, then $A^\#$ exists.
\end{lemma}

\begin{beweis}
   Let $X \prec (H_\theta;\in)$ where $\theta >> \aleph_{\omega + \omega}$, $A \in X$, $\card(X) = \aleph_1$, $X$ is countably closed, and $X$ is cofinal in $(\card(A)^+)^{L(A)}$. Let $\sigma_X:H_X \rightarrow X$ be an isomorphism where $H_X$ is transitive. Let $\bar{A} = \sigma^{-1}_X(A)$. $\bar{A}$ is coded by a subset of $\aleph_1$, so $\bar{A}^\#$ exists. A standard lift-up argument using countable closure of $X$ then gives the existence of $A^\#$.
\end{beweis}

\begin{lemma}
    Let $A \subseteq \aleph_\omega$. Then $A^\#$ exists.
\end{lemma}

\begin{beweis}
   Let $X \prec (H_\theta;\in)$ where $\theta >> \aleph_{\omega + \omega}$, $A \in X$, and $X \cap H_{\aleph_\omega}$ is good. Let $\sigma_X: H_X \rightarrow X$ be an isomorphism where $H_X$ is transitive. We can assume that $\cof((\card(A)^+)^{L(A)}) \geq \omega_1$. Thus by elementarity and the closure of $X$ we have \[\cof((\card(\bar{A})^+)^{L(\bar{A})}) \geq \aleph_1\] where $\bar{A} := \sigma^{-1}_X(A)$. Note that we do not need to write $(L(\bar{A})^{H_X})$ here as $\card(\bar{A})^+ \subset H_X$. Now $\bar{A}$ has size $\aleph_1$ so $\bar{A}^\#$ exists. But this implies $\cof((\card(\bar{A})^+)^{L(\bar{A})}) = \omega$. Contradiction!
\end{beweis}

Using this lemma and applying the same argument again we can immediately conclude:

\begin{lemma}
   Let $A \subset \aleph_{\omega + \omega}$ be bounded. Then $A^\#$ exists.
\end{lemma}

This finishes the first step of the induction. Assume then that we already know that $M^\#_n(A)$ exists for all bounded subsets $A \subset \aleph_{\omega + \omega}$. For the argument we need our own version of the core model dichotomy. Let $A \subset \aleph_{\omega + \omega}$ be a bounded subset. We let $K^c(A)$ be the result of the robust extender construction above $A$ using extenders with critical point below $\aleph_{\omega + \omega}$ if it exists, i.e. all levels of the construction have solid standard parameters.

\begin{lemma}
    Let $A \subset \aleph_{\omega + \omega}$ be a bounded subset. Then (exactly) one of the following is true:
    \begin{itemize}
        \item[$(i)$] $K^c(A)$ exists, is $(n+1)$-small and is $(\aleph_{\omega + \omega})$-iterable;
        \item[$(ii)$] some level of the robust extender construction is not $(n+1)$-small and is $(\aleph_{\omega + \omega})$-iterable.
    \end{itemize}
\end{lemma}

\begin{beweis}
   Let $\mathcal{M}$ be a level of the robust extender construction over $A$ that is $(n+1)$-small. We want to see that $\mathcal{M}$ is countably iterable. To that end consider $\mathcal{M}^+ := M^\#_n(\mathcal{M})$. Any countable hull of $\mathcal{M}^+$ is then uniquely iterable for countable trees using only extenders from $\mathcal{M}$ by the realizable branches strategy. A reflection argument like \thref{strrefl} then shows that $\mathcal{M}^+$ is $(\aleph_{\omega + \omega})$-iterable for trees using only extenders from $\mathcal{M}$. This implies that $\mathcal{M}$ is sufficiently iterable to have solid parameters. Thus if all levels of the construction are $(n+1)$-small, then $(i)$ holds.

Assume then $\mathcal{M}$ is the least level of the construction that is not $(n+1)$-small. Then any countable hull of $\mathcal{M}$ is uniquely iterable for countable trees using the realizable branches strategy which is thus the same as the $M^\#_n$-guided iteration strategy. \thref{strrefl} then gives that $\mathcal{M}$ is uniquely $(\aleph_{\omega + \omega})$-iterable, so $(ii)$ holds.
\end{beweis}

Let now $A \subseteq \aleph_1$. Our next goal is to show that $M^\#_{n+1}(A)$ exists. Assume not. Then by the core model dichotomy $K^c(A)$ exists and is $(\aleph_{\omega + \omega})$-iterable. We will now construct not one, but a countable sequence of soundness witnesses. Let $\Omega_m := \aleph_{\omega + (m + 5)}$. An $m$-weasel is a weasel of height $\Omega_m$. For an $m$-weasel $W$ let $W^+ := S(W)$ be the stack over $W$. Note that $K^c(A)\vert\vert \Omega_m$ is in the language of \cite{knombl} a ``mini-universe" by \thref{nosquarelemma} and \cite{square}. Therefore we can extract from $(K^c(A)\vert\vert\Omega_m)^+$ an $m$-weasel $W_m$ such that $W^+_m$ has the definability property for all $\alpha < \aleph_{\omega + (m + 3)}$. 

Let $X \prec (H_\theta;\in)$ be such that $A \in X$, $\<W_m: m < \omega\> \in X$, and $X \cap H_{\omega + \omega}$ is good. Let $\sigma_X: H_X \rightarrow X$ be an isomorphism where $H_X$ is transitive. We let $W^X_m := \sigma^{-1}_X(W_m)$. We will form a sequence of trees: $\mathcal{T}_m$ is the tree on $W_m$ arising from the co-iteration with $W^X_m$ but only up to $\aleph_{m + 3}$.

Let $\mathcal{M}^m_\beta$ be the $\beta$-th model in $\mathcal{T}_m$, $E^m_\beta$ the $\beta$-th extender, $\iota^m_{\beta,\gamma}$ the iteration embedding, $\mathcal{D}_m$ the set of drops. Note that by \cite[Section 6]{knombl} we have $W_m \vert\vert\aleph_{\omega + (m + 3)} = W_k \vert\vert\aleph_{\omega + (m + 3)}$ for all $m \leq k < \omega$. This implies by an easy induction:

\begin{lemma}\label{sec4agreement}
    $\leq_{\mathcal{T}_m} = \leq_{\mathcal{T}_k} \restr \length(\mathcal{T}_m)$, $E^m_\beta = E^k_\beta$ for all $\beta < \length(\mathcal{T}_m)$, and $\mathcal{D}_m = \mathcal{D}_k \cap \length(\mathcal{T}_m)$ for all $m \leq k < \omega$.
\end{lemma}

Moreover each $W_m$ satisfies the weak covering property for each $\alpha \in \left[\aleph_2,\aleph_{\omega + (m + 3)}\right)$.This implies just as in \thref{sec3cofincore} that:

\begin{lemma}\label{sec4cofincore}
   Let $\kappa < \aleph_{m + 3}$ be above $\crit(\sigma_X)$ for some $m <  \omega$. Then $\cof((\kappa^+)^{W^X_m}) = \card(\kappa)$.
\end{lemma}

Let now $\kappa_\alpha := \aleph^{W^X_m}_\alpha$ and $\lambda_\alpha := \aleph^{W^X_m}_{\alpha + 1}$ where $m$ is such that $\alpha < \aleph_{m+3}$. We can then define models $\mathcal{P}^m_\alpha,\mathcal{Q}^m_\alpha,\mathcal{R}^m_\alpha,\mathcal{S}^m_\alpha$ for $\alpha < \aleph_{m+3}$ derived from the tree $\mathcal{T}_m$ as in the last section.  In general we will have $\mathcal{P}^m_\alpha \neq \mathcal{P}^k_\alpha$ etc, but the structure of the phalanxes will be compatible, i.e. if $\mathcal{P}^m_\alpha \eextend \mathcal{M}^m_\beta$ then $\mathcal{P}^k_\alpha \eextend \mathcal{M}^k_\beta$ for all $m \leq k <\omega$ and $\alpha < \aleph_{m + 3}$. We will have the following properties for $m < \omega$ and all $\alpha < \aleph_{m + 3}$:

\begin{itemize}
    \item $(1)^m_\alpha$ is the statement that the least extender applied to $W^X_m$ in the co-iteration of $W^X_m$ and $W_m$ has length at least $\lambda_\alpha$;
    \item $(2)^m_\alpha$ is the statement that the phalanx $\<\<W_m,\mathcal{S}^m_\alpha\>,\sigma_X(\kappa_\alpha)\>$ is iterable;
    \item $(3)^m_\alpha$ is the statement that the phalanx $\<\<W^X_m,\mathcal{Q}^m_\alpha\>,\kappa^X_\alpha\>$ is iterable;
    \item $(4)^m_\alpha$ is the statement that the phalanx $\<\<\mathcal{P}^m_\beta: \beta < \alpha\> \concat W^X_m\>,\<\lambda_\beta: \beta < \alpha\>\>$ is iterable;
    \item $(5)^m_\alpha$ is the statement that the special phalanx $\<\<\mathcal{R}^m_\beta: \beta < \alpha\> \concat W_m\>,\<\chi_X(\sigma_X(\lambda_\beta)):\beta < \alpha\>\>$ is iterable;
    \item $(6)^m_\alpha$ is the statement that the special phalanx $\<\<\mathcal{S}^m_\beta: \beta < \alpha\> \concat W_m\>,\<\chi_X(\sigma_X(\lambda_\beta)):\beta < \alpha\>\>$ is iterable.
\end{itemize}

This is proved just like as in \cite{coveringbelowwoodin}. Note that just as in \cite[4.22]{knombl} the iterability of the phalanx $\<\<W^X_m,\mathcal{S}^X_\alpha\>,\sigma_X(\kappa_\alpha)\>$ implies the iterability of $\<\<(W^X_m)^+,(\mathcal{S}^X_\alpha)^+\>,\sigma_X(\kappa_\alpha)\>$ so we can make use of definability property in the usual fashion.

Let then $\alpha_k$ be the unique $\alpha$ such that $\kappa_{\alpha_k} = \aleph_{k+2}$. Let $\xi_n$ be least $\xi$ such that $\nu(E^m_\xi) \geq \aleph_{n + 2}$ for $m \geq k$ (if such exists otherwise $\xi_k = \length(\mathcal{T}_m)$). The next then follows just as in \thref{sec3cutpoint}:

\begin{lemma}\label{sec4cutpoint}
    $\beta \geq \xi_k \Rightarrow \beta \geq_{\mathcal{T}_m} \xi_k$ for all $k \leq m < \omega$.
\end{lemma}

Next we get our first drop just as in \thref{sec3firstdrop}

\begin{lemma}
    $\mathcal{P}^m_{\alpha_k}$ is not a weasel for all $k \leq m < \omega$.
\end{lemma}

Note that $\mathcal{P}^m_{\alpha_0}$ is a mouse of size $\aleph_2$ so by agreement between the iteration trees we have $\mathcal{P}^m_{\alpha_0} = \mathcal{P}^k_{\alpha_0}$ for $m \leq k < \omega$. Call this mouse $\mathcal{N}$. By \thref{sec4agreement} $\mathcal{T}_m$ can then be considered as the concatenation $\mathcal{T}^*_m\concat \mathcal{U}_m$ where $\mathcal{T}_m$, a tree on $W_m$, only varies in models and a tree $\mathcal{U}_m$ on $\mathcal{N}$. For $m \leq k < \omega$ the tree $\mathcal{U}_k$ is an extension of the tree $\mathcal{U}_m$. So we can form $\mathcal{U} = \union_{m < \omega} \mathcal{U}_m$. $\mathcal{U}$ has a cofinal branch determined by $\{\xi_k :  k < \omega\}$. We are done if we can show that there are infinitely many drops on this branch as this contradicts the iterability of $\mathcal{N}$. But this exact thing follows as in \thref{sec3moredrops}. 

We shall now move on to subsets $A \subseteq \aleph_\omega$. Our goal is to prove an instance of mouse reflection, i.e. we are going to leverage the just proven existence of $M^\#_{n + 1}$'s for subsets of $\aleph_1$. Yet the argument to come follows the same general outline.

Let us assume the that $M^\#_{n + 1}(A)$ does not exist by which we mean that no level of the robust $K^c(A)$-construction is not $(n+1)$-small. We can thus define the sequence $\<W_m : m < \omega\>$ as before. These models satisfy a version of weak covering, i.e. $\cof((\alpha^+)^{W_m}) \geq \card(\alpha)$ where $\alpha \in \left[\aleph_{\omega + 1},\aleph_{\omega + (m + 3)}\right)$ and $m < \omega$. This is seen by the same method from \cite{coveringbelowwoodin} with some obvious changes.

Let $X \prec (H_\theta;\in)$ be such that $A \in X$, $\<W_m: m < \omega\> \in X$, and $X \cap H_{\omega + \omega}$ is good. Let $\sigma_X: H_X \rightarrow X$ be an isomorphism where $H_X$ is transitive. We let $A_X := \sigma^{-1}_X(A)$ and $W^X_m := \sigma^{-1}_X(W_m)$. This is as before, but it is here where the arguments must diverge. We cannot co-iterate $W^X_m$ with $W_m$ as in the usual covering set up as they do not belong to the same hierarchy of mice. Note though that $\card(A_X) = \aleph_1$. Hence $M^\#_{n+1}(A_X)$ exists and is not $(n+1)$-small. Therefore by assumption it must win the co-iteration with $W^X_m$. The idea then is to (partially) substitute $M^\#_{n+1}(A_X)$ for $W_m$ in the covering argument. We will have to see that this works.

Let $\mathcal{T}_m$ be the tree on $M^\#_{n+1}(A_X)$ that arises from the co-iteration with $W^X_m$ but only up to $\aleph_{m + 3}$. We will use the previously established notation for models, extenders, etc. We also have that $W_m \vert\vert \aleph_{\omega + (m + 3)} = W_k \vert\vert \aleph_{\omega + (m + 3)}$ for all $m \leq k < \omega$. This easily implies:

\begin{lemma}
    $\mathcal{T}_k$ end-extends $\mathcal{T}_m$ for all $m \leq k < \omega$.
\end{lemma}

This is actually slightly better than \thref{sec4agreement} as $\mathcal{T}_m$ are based on the same model. Our eventual contradiction will come from the fact that $\union_{m < \omega} \mathcal{T}_m$, a tree on $M^\#_{n+1}(A_X)$, has a canonical cofinal branch with infinitely many drops on it. We do so by using the weak covering of the $W_m$. 

\begin{lemma}\label{sec4cofincore2}
   Let $\kappa < \aleph_{m + 3}$ be above $\aleph_2$ for some $m <  \omega$. Then $\cof((\kappa^+)^{W^X_m}) = \card(\kappa)$.
\end{lemma}

Unfortunately, this only applies to $W^X_m$ and not its iterates, so we will have to maintain that $W^X_m$ does not move in the co-iteration. The only way we can see to prove this is by going through the covering argument with its 6 properties even though we only really care about the first one. 

Let $\alpha < \aleph_{m + 3}$ for any $m$: let $\kappa_\alpha := \aleph^{W^X_m}_\alpha$ and $\lambda_\alpha := \aleph^{W^X_m}_{\alpha + 1}$. Let $\mathcal{P}^m_\alpha$ be the least initial segment of $\mathcal{M}^m_\beta$ that agrees with $W^X_m$ up to $\lambda_\alpha$ and projects to or below $\kappa_\alpha$ where $\beta$ is minimal with $\length(E^m_\beta) \geq \lambda_\alpha$. $n^m_\alpha$ is least such that $\rho_{n_\alpha + 1}(\mathcal{P}^m_\alpha) \leq \kappa_\alpha$. Let $\mathcal{R}^m_\alpha := \ult_{n^m_\alpha}(\mathcal{P}^m_\alpha; \sigma_X \restr (W^X_m \vert\vert\lambda_\alpha))$. If $\mathcal{R}^m_\alpha$ is a proto-mouse we let $\mathcal{S}^m_\alpha := \ult_{k^m_\alpha}(\mathcal{S}^m_\beta;F^{\mathcal{R}^m_\alpha})$ where $\beta = \mu^{\mathcal{R}^m_\alpha}$ and $k^m_\alpha = k^m_\beta$, otherwise $\mathcal{S}^m_\alpha := \mathcal{R}^m_\alpha$ and $k^m_\alpha = n^m_\alpha$. If $\mathcal{R}^m_\alpha$ is a protomouse we will also let $\mathcal{Q}^m_\alpha := \ult_{k^m_\alpha}(\mathcal{Q}^m_\beta; F^{\mathcal{P}^m_\alpha})$, otherwise $\mathcal{Q}^m_\alpha := \mathcal{P}^m_\alpha$.

It is not hard to see that $\mathcal{P}^m_\alpha = \mathcal{P}^k_\alpha$ if $m \leq k < \omega$ and $\alpha < \aleph_{m + 3}$. Never the less we will keep the superscript as each $m$ will have its own induction along the $\alpha$.

The six properties which should look familiar at this point are then:

\begin{itemize}
    \item $(1)^m_\alpha$ is the statement that the least extender applied to $W^X_m$ in the co-iteration of $W^X_m$ and $W_m$ has length at least $\lambda_\alpha$;
    \item $(2)^m_\alpha$ is the statement that the phalanx $\<\<W_m,\mathcal{S}^m_\alpha\>,\sigma_X(\kappa_\alpha)\>$ is iterable;
    \item $(3)^m_\alpha$ is the statement that the phalanx $\<\<W^X_m,\mathcal{Q}^m_\alpha\>,\kappa^X_\alpha\>$ is iterable;
    \item $(4)^m_\alpha$ is the statement that the phalanx $\<\<\mathcal{P}^m_\beta: \beta < \alpha\> \concat W^X_m\>,\<\lambda_\beta: \beta < \alpha\>\>$ is iterable;
    \item $(5)^m_\alpha$ is the statement that the special phalanx $\<\<\mathcal{R}^m_\beta: \beta < \alpha\> \concat W_m\>,\<\chi_X(\sigma_X(\lambda_\beta)):\beta < \alpha\>\>$ is iterable;
    \item $(6)^m_\alpha$ is the statement that the special phalanx $\<\<\mathcal{S}^m_\beta: \beta < \alpha\> \concat W_m\>,\<\chi_X(\sigma_X(\lambda_\beta)):\beta < \alpha\>\>$ is iterable.
\end{itemize}

The covering argument mainly uses properties of the core model: the definability property, and the iterability of core model derived phalanxes. We will be able to use the same argument then for $M^\#_{n+1}(A_X))$ by using soundness, and the iterability of $M^\#_{n+1}(A_X)$ derived phalanxes. (We will have that all $\mathcal{S}^m_\alpha$ are mice, so we have no use for the hull property.)

\begin{definition}
  A phalanx $\<\<\mathcal{N}_\alpha: \alpha < \beta + 1\>,\<\gamma_\alpha: \alpha < \beta\>\>$ is $M^\#_{n+1}(A_X)$ derived iff there exists a sequence $\<\mathcal{U}_\alpha : \alpha < \beta + 1\>$ of iteration trees such that for all $\alpha < \delta < \beta + 1$ there is some $\xi$ such that $\mathcal{U}_\alpha \restr \xi = \mathcal{U}_\delta \restr \xi$ and $\gamma_\alpha \leq \min(\length(E^{\mathcal{U}_\alpha}_\xi),\length(E^{\mathcal{U}_\delta}_\xi))$; if one, or both, of these extenders is not defined then substitute $\on$.
\end{definition}

\begin{lemma}
   Any $M^\#_{n+1}(A_X)$ derived phalanx is iterable.
\end{lemma}

For this lemma it is important to recall how we showed the existence of the $M^\#_{n+1}(A_X)$ in the first place by the use of the Core Model Dichotomy, i.e. there is some level $\mathcal{M}$ of the robust $K^c(A_X)$ construction such that $\mathcal{C}_\omega(\mathcal{M}) = M^\#_{n+1}(A_X)$. Let us write $\pi$ for the reverse core embedding.

\begin{beweis}
   Assume not. Take some countable $Y \prec (H_\theta;\in)$ containing $\mathcal{I}$ a counter example. Let $\sigma_Y:H_Y \rightarrow Y$ be an isomorphism where $H_Y$ is transitive. Write $M_Y := \sigma^{-1}_Y(M^\#_{n+1}(A_X))$ and $\mathcal{I}_Y := \sigma^{-1}_Y(\mathcal{I})$.
   
   By a \thref{strrefl}-like reflection argument the phalanx $\mathcal{I}_Y$ fails to be iterable. On the other hand by elementarity it is $M_Y$-derived. $M_Y$ is embeddable into $\mathcal{M}$ by $\pi \circ \sigma_Y$. Thus by \cite{robust} $\mathcal{I}_Y$ can be lifted onto $\mathcal{M}$. Contradiction!
\end{beweis}

These are all the ingredients we will need apart from the closure of $X$. Let us now fix some $m < \omega$ and some $\alpha < \aleph_{m + 3}$. Assume that $(1)^m_\beta$ through $(6)^m_\beta$ hold for all $\beta < \alpha$.

\begin{lemma}
    For all $\beta < \alpha$ there exists an iteration tree $\mathcal{U}_\beta$ on $W_m$ and some $r\Sigma_{n_\beta + 1}$- or $r\Sigma_{k_\alpha + 1}$-elementary embedding (depending on how $\mathcal{S}^m_\beta$ was constructed) $\upsilon_\beta$ from $\mathcal{S}^m_\beta$ into the last model of $\mathcal{U}_\beta$ or some initial segment thereof. We have that $\length(E^{\mathcal{U}_\beta}_0) \geq \chi_X(\sigma_X(\lambda_\beta))$ and $\crit(\upsilon_\beta) \geq \sigma_X(\kappa_\beta)$.
\end{lemma}

Note we should actually have here that $\mathcal{S}^m_\beta$ is an initial segment of $W_m$ but strictly speaking we do not need it.

\begin{beweis}
    Fix $\beta < \alpha$. By $(2)^m_\beta$ the phalanx $\<\<W_m,\mathcal{S}^m_\beta\>,\sigma_X(\kappa_\beta)\>$ is iterable. We can co-iterate it against $W_m$. By a standard argument the final model of the iteration on the phalanx side is above $\mathcal{S}^m_\beta$. Thus the resulting iteration embedding which we call $\upsilon_\beta$ has critical point at least $\sigma_X(\kappa_\beta)$. Finally note that $W_m$ and $\mathcal{S}^m_\beta$ agree up to $\chi_X(\sigma_X(\lambda_\beta))$. 
\end{beweis}

The main corollary here is that the phalanx $\<\<\mathcal{M}^{\mathcal{U}_\beta}_{\length(\mathcal{U}_\beta) - 1}: \beta < \alpha\> \concat \<W_m\>,\<\chi_X(\sigma_X(\lambda_\beta)):\beta < \alpha\>\>$ is $W_m$-derived and hence iterable. The phalanx of $(6)^m_\alpha$ can then be lifted onto this phalanx via the embeddings $\<\upsilon_\beta:\beta < \alpha\>$.

$(5)^m_\alpha$ can then be derived the same way from $(6)^m_\alpha$ as in \cite[3.18]{coveringbelowwoodin}. Nothing about this particular argument is specific to the core model so we shall skip further detail.

$(4)^m_\alpha$ then follows immediately as the phalanx can be copied onto the phalanx of $(5)^m_\alpha$ via the assorted lift-up maps.

\begin{lemma}
   For all $\beta \leq \alpha$ there exists an iteration tree $\mathcal{V}_\beta$ on $M^\#_{n+1}(A_X)$ that end-extends $\mathcal{T}_m \restr \gamma_\beta$ where $\gamma_\beta$ is least such that $\mathcal{P}^m_\beta \eextend \mathcal{M}^m_{\gamma_\beta}$. We have that $\length(E^{\mathcal{V}_\beta}_{\gamma_\beta}) \geq \lambda_\beta$. Furthermore there exists an embedding $\psi_\beta$ from $W^X_m$ into the final model of $\mathcal{V}_\beta$ or an initial segment thereof. $\crit(\psi_\beta) \geq \kappa_\beta$.
\end{lemma}

\begin{beweis}
   We will use $(4)^m_\beta$. We co-iterate the phalanxes $\<\<\mathcal{P}^m_\delta: \delta < \beta\>\concat\<W^X_m\>,\<\lambda_\delta:\delta < \beta\>\>$ and $\<\<\mathcal{P}^m_\delta : \delta \leq \beta\>,\<\lambda_\delta:\delta < \beta\>$. Note that the first disagreement between $\mathcal{P}^m_\beta$ and $W^X_m$ by design is at least $\lambda_\beta$. So the tree on the latter phalanx is our desired $\mathcal{V}_\beta$. If we can show that the iteration on the first phalanx is above $W^X_m$ then we are done as the iteration embedding will be our $\psi_\beta$.
   
 So assume that this is not the case, say the final model is above $\mathcal{P}^m_\delta$ for some $\delta < \beta$. Note that this implies that some extender is applied during the iteration to it as otherwise we are above $W^X_m$ by default. This implies that the final model on the first phalanx side, call it $\mathcal{M}$, is not sound. By standard arguments the iteration embedding from the last drop into the last models is simply the reversed core embedding.
 
 This is important as this implies that no extender can be applied on the second phalanx. Otherwise the same argument applies and we must conclude that compatible extenders are applied in the course of the iteration. This is a contradiction as usual for comparisons.
 
 We then have that $\mathcal{P}^m_\beta \eextend \mathcal{M}$, but in fact equality must hold. The key fact is that by choice $\mathcal{P}^m_\beta$ defines a subset of $\kappa_\beta$ that is not in $W^X_m$ and thus also not in $\mathcal{M}$.
 
 We will compare two branch tails: the first, $t_0$ is the tail of the branch that leads up to $\mathcal{P}^m_\delta$ in $\mathcal{T}_m$ concatenated with the branch that leads from $\mathcal{P}^m_\delta$ up to $\mathcal{M}$ from its last drop onwards; the second, $t_1$, is the branch of $\mathcal{T}_m$ leading up to $\mathcal{P}^m_\beta$ from its last drop onwards; note that the former sequence may include extenders from $\mathcal{T}_m$ too. We want to show that both of these tails agree on all models and branch embeddings, if not their indexes.
 
 Both tails lead up to $\mathcal{M} = \mathcal{P}^m_\beta$ starting from their core. Assume not: let $E,F$ be the first two extenders, $E$ used on $t_0$ and $F$ on $t_1$, that are not equal; by their minimality both extenders are applied to the same model ; in fact, both equal the core above the supremum of generators of the previously applied extenders and thus their cumulative embeddings are equal; we must therefore have that $E,F$ are compatible; as we assumed that they are not equal one is a proper initial segment of the other, say $E$ is an initial segment of $F$; by the initial segment condition therefore $E \in \mathcal{M}$, but on the other hand $\length(E)$ is a cardinal there; contradiction!
 
 But this is absurd, as the first extender used during the iteration is longer than any extender that was used in $\mathcal{T}_m \restr \gamma_\beta$.
\end{beweis}

We can now prove $(1)^m_\alpha$. Notice that the existence of $\psi_\alpha$ implies that $\Pot^{W^X_m}(\kappa_\alpha)$ is contained in the final model of $\mathcal{V}_\alpha$. Therefore the least disagreement between $W^X_m$ and $\mathcal{P}^m_\alpha$ is above $\lambda_\alpha$. By the choice of $\mathcal{P}^m_\alpha$ the same is true for the least disagreement between $\mathcal{M}^m_{\gamma_\alpha}$ and $W^X_m$.

Even though we were really only interested in proving $(1)^m_\alpha$ we still have to close the loop to finish the induction. Thankfully we can get $(3)^m_\alpha$ easily now. The phalanx $\<\<\mathcal{M}^{\mathcal{V}_\beta}_{\length(\mathcal{V}_\beta)-1}\>\concat\<\mathcal{Q}^m_\alpha\>,\<\kappa_\beta:\beta < \alpha\>$ is $\mathcal{M}^\#_{n+1}(A_X)$-derived and hence iterable. We can then copy a tree on the phalanx $\<\<W^X_m,\mathcal{Q}^m_\alpha\>,\kappa_\alpha\>\>$ onto a tree on this phalanx using the maps $\<\psi_\beta: \beta < \alpha\>\concat\id$. The way this works is non-standard but quite simple: any extender applied to $W^X_m$ in the  course of the iteration must have critical point some $\kappa_\beta$ for $\beta < \alpha$, this extender is then copied to an extender over $\mathcal{M}^{\mathcal{V}_\beta}_{\length(\mathcal{V}_\beta)-1}$ using the map $\psi_\beta$, the rest is as usual.

$(2)^m_\alpha$ follows then by a standard argument using the countable closure of $X$ just like in \cite[3.13]{coveringbelowwoodin}. This completes the induction step.

We thus have $(1)^m_\alpha$ for all $\alpha < \aleph_{m + 3}$ and $m < \omega$. We can now proceed to the main argument which is in fact just a simplification of previous arguments. We will form the tree $\mathcal{T} := \union_{m < \omega} \mathcal{T}_m$ on $M^\#_{n+1}(A_X)$. This tree has a canonical cofinal branch. In fact, no extender $E$ can be applied during this iteration tree with the property that $\crit(E) < \aleph_{m + 2}$ but $\length(E) \geq \aleph_{m + 2}$, any $m < \omega$. Otherwise $\length(E)$ is a successor cardinal in $W^X_m$ and hence by \thref{sec4cofincore2} has cofinality at least $\aleph_{m + 2}$, on the other hand that cofinality equals the cofinality of $\crit(E)^+$ which is below $\aleph_{m + 2}$.

Notice then that we start with a model of size $\aleph_1$ so for any $m < \omega$ we have \[\cof((\aleph^+_{m + 2})^{\mathcal{M}^m_{\aleph_{m + 2}}}) < \aleph_{m + 2}.\] But by \thref{sec4cofincore2} $\cof((\aleph^+_{m + 2})^{W^X_m}) \geq \aleph_{m + 2}$. Hence any extender applied to $\mathcal{M}^m_{\aleph_{m + 2}}$ in the course of the iteration induces a drop. But by the previous fact the cofinal branch through $\mathcal{T}$ passes through all these nodes so contains infinitely many drops. Contradiction!

So we can conlude that $M^\#_{n + 1}(A)$ does exist. We are finished if we can show the same thing for $A$ that is a bounded subset of $\aleph_{\omega + \omega}$. The argument is mostly the same. Given such $A$ we form soundness witnesses $W_m$. Depending on the exact size of $A$ these will only be defined for all but finitely many $m$, but this makes no difference. We then take some appropriate $X$ containing $A$. $A$ will collapse down to some $A_X$ of size less than $\aleph_\omega$. By what we know already $M^\#_{n + 1}(A_X)$ exists. We then co-iterate it with the collapses of the $W_m$'s. This leads to the same contradiction as before. We shall skip further detail, finishing the induction step and hence the proof of $\PD$.

\section{Open Questions}

An important question is if similar strength can be derived from a simpler version of Chang's Conjecture such as $(\aleph_3,\aleph_2) \twoheadrightarrow (\aleph_2,\aleph_1)$. Unfortunately, our argument for an inner model with a Woodin cardinal is quite specific to the ``long" version of the Chang's Conjecture.

\begin{frage}\label{fragedrei}
  Assume $(\aleph_3,\aleph_2) \twoheadrightarrow (\aleph_2,\aleph_1)$ and $\GCH$. Does there exist an inner model with a Woodin cardinal?
\end{frage}

We think though we might know how the argument would proceed from an inner model with a Woodin cardinal. We think the core model induction ought to be done in $V^{\col(\omega,\omega_2)}$. Consider, say, $(\mathcal{P},\Sigma)$ a $\Gamma$-suitable pair for some determined pointclass $\Gamma$ from $V^{\col(\omega,\omega_2)}$. This pulls back to some pair, which we call by the same name, over $V$ where $\mathcal{P}$ has size $\kleinergleich \aleph_2$ and $\Sigma$ is a $\aleph_3$-iteration strategy. The first step would be to show that for every hull $\sigma: \bar{\mathcal{P}} \rightarrow \mathcal{P}$ of size less than $\aleph_2$ we have that $M^{\Sigma^\sigma,\#}_n(A)$ exists for all $n$ and $A$ a subset of $\aleph_2$. This would presumably be a straightforward adaption of whatever argument yields $\PD$. Then if some $M^{\Sigma,\#}_n(A)$ fails to exists, then there must be some stable core model $K^\Sigma(A)$ of height $\aleph_3$. Consider some Chang structure $X$ containg this model. It's collapse $K^\Sigma_X(A_X)$ could then be compared with $M^{\Sigma^{\sigma_X},\#}_n(A_X)$, hopefully yielding a contradiction in the same manner.

Could we apply the same idea to the property covered in this work? Unfortunately, this seems to be difficult. We should want to do the core model induction in a symmetric extension of the universe where $\aleph^V_{\omega + \omega}$ is $\omega_1$. But in such a model the principle of dependent choice fails. This principle is needed for the core model induction method as we understand it today.

One could try to operate similar to \cite{SteelZoble} and show that given some $\alpha_n$ beginning a gap in $V^{\col(\omega,\aleph_{\omega + n})}$ there must be $\alpha_m$ begininning a corresponding gap in $V^{\col(\omega,\aleph_{\omega + m})}$ for some $n \leq m < \omega$. This would likely yield success.

\begin{frage}
	Assume $(\aleph_{\omega + 1},\aleph_{\omega + 2},\ldots) \twoheadrightarrow (\aleph_1,\aleph_2,\ldots)$ and $\GCH$. Does $\AD^{L(\mathbb{R})}$ hold?
\end{frage}

In \cite{OpistolfromJonsson} we were able to prove the existence of $0^\P$ from the weaker assumption of $(\aleph_{\omega + \omega},\aleph_\omega) \twoheadrightarrow (\aleph_\omega, \aleph_n)$ for some $n < \omega$. Once we move past the level of linear iterations it seems we need extra assumptions to make sure that we have infinitely many drops occuring on some branch. It is unclear to us if these extra assumptions can be removed.

\begin{frage}
  Assume $(\aleph_{\omega + \omega},\aleph_\omega) \twoheadrightarrow (\aleph_\omega, \aleph_n)$ for some $n < \omega$. Does there exist an inner model with a Woodin cardinal?
\end{frage}

\bibliographystyle{alpha}
\bibliography{bibliographie}

\end{document}